\documentclass[12pt, english]{amsart}

\usepackage{geometry}
\usepackage[english]{babel}
\usepackage{multirow}% pacote para matrizes em blocos
\usepackage{graphicx}
\usepackage{amssymb}
\usepackage{yhmath}
\usepackage{verbatim} 
\usepackage{amsmath}
\usepackage{amsthm}
\usepackage{color}
\usepackage{epstopdf} %Figuras .eps
\usepackage{wrapfig}
\usepackage[all]{xy}
\usepackage{hyperref}%links internos
\usepackage{enumitem} %para letras e outros símbolos no ambiente 'enumerate'

\newcommand{\nn}{\nonumber}
\newcommand{\vol}{{\rm vol}}

\newtheorem{theorem}{Theorem}
\newtheorem{corollary}{Corollary}

\newtheorem{definition}{Definition}
\newtheorem*{theoremA}{Theorem A}
\newtheorem*{theoremB}{Theorem B}

 \theoremstyle{remark}

\begin{document}
\title[Normal degree and a lower bound for unit vector fields on hypersurfaces]{A topological lower bound for the energy of a unit
 vector field on a closed Euclidean hypersurface}
\author{Fabiano G. B. Brito}

\author{Icaro Gon\c {c}alves}

\author{Adriana V. Nicoli}

\address{Centro de Matem\'atica, Computa\c{c}\~ao e Cogni\c{c}\~ao,
Universidade Federal do ABC,
09.210-170 Santo Andr\'e, Brazil}
\email{fabiano.brito@ufabc.edu.br}

\address{Dpto. de Matem\'{a}tica, Instituto de Matem\'{a}tica e Estat\'{i}stica,
Universidade de S\={a}o Paulo, R. do Mat\={a}o 1010, S\={a}o Paulo-SP
05508-900, Brazil.}
\email{avnicoli@ime.usp.br}

\address{Dpto. de Matem\'{a}tica, Instituto de Matem\'{a}tica e Estat\'{i}stica,
Universidade de S\={a}o Paulo, R. do Mat\={a}o 1010, S\={a}o Paulo-SP
05508-900, Brazil.}
\email{icarog@ime.usp.br}

%57R25 Vector fields, frame fields
%47H11 Degree Theory
%58E20 Harmonic maps

\thanks{Icaro Gon\c calves was supported by CNPq 141113/2013-8, and is partially supported by a scholarship from the National Postdoctoral Program, PNPD-CAPES. Adriana V. Nicoli is supported by CAPES-PROEX 1745551}

\begin{abstract}
For a unit vector field on a closed immersed Euclidean hypersurface $M^{2n+1}$, $n\geq 1$, we exhibit a nontrivial lower bound for its energy which depends on the degree of the Gauss map of the immersion. When the hypersurface is the unit sphere $\mathbb{S}^{2n+1}$, immersed with degree one, this lower bound corresponds to a well established value from the literature. We introduce a list of functionals $\mathcal{B}_k$ on a compact Riemannian manifold $M^{m}$, $1\leq k\leq m$, and show that, when the underlying manifold is a closed hypersurface, these functionals possess similar properties regarding the degree of the immersion. In addition, we prove that Hopf flows minimize $\mathcal{B}_n$ on  $\mathbb{S}^{2n+1}$.
\end{abstract}
\maketitle

\section{Introduction and statement of the main results}
Let $M^m$ be a compact oriented Riemannian manifold, $m\geq 2$, and let $\nabla$ denote its Levi-Civita connection. The energy of a unit vector field on $M$ is defined as the energy of the map $\vec{v}:M\rightarrow T_1M$, where $T_1M$ denotes the unit tangent bundle equipped with the Sasaki metric, (see \cite{wiegmink} and \cite{wood})
\begin{eqnarray}
\label{energy1}
E(\vec{v})=\frac{1}{2}\int_M\|\nabla\vec{v}\|^2+\frac{m}{2}\vol(M).
\end{eqnarray}
In \cite{wiegmink}, Wiegmink defines the total bending functional, a quantitative measure for the extent to which a unit vector field fails to be parallel with respect to the Levi-Civita connection  $\nabla$ of a Riemannian manifold $M$. Precisely, 
\begin{eqnarray}
\mathcal{B}(\vec{v})=\frac{1}{(m-1)\vol(\mathbb{S}^m)}\int_M\|\nabla\vec{v}\|^2, \nn
\end{eqnarray}
and the energy of $\vec{v}$ may be written in terms of this functional as
\begin{eqnarray}
E(\vec{v})=\frac{(m-1)\vol(\mathbb{S}^m)}{2}\mathcal{B}(\vec{v})+\frac{m}{2}\vol(M).\nn
\end{eqnarray}

An important question regarding these functionals is whether one can find unit vector fields such that the minimum of the above functional is attained. Brito \cite{brito} showed that Hopf flows are absolute minima of the functional $\mathcal{B}$ in $\mathbb{S}^3$:
\begin{theorem}[Brito, \cite{brito}]
Hopf vector fields are the unique vector fields on $\mathbb{S}^3$ to minimize $\mathcal{B}$.
\end{theorem}
Gluck and Ziller proved that Hopf flows are also the unit vector fields of minimum volume, with respect to the following definition of volume,
$$\vol(\vec{v})=\int_M\sqrt{\det(I+(\nabla\vec{v})(\nabla\vec{v})^*)},$$
where $I$ is the identity and $(\nabla\vec{v})^*$ represents adjoint operator. 
\begin{theorem}[Gluck and Ziller, \cite{GZ}]
The unit vector fields of minimum volume on $\mathbb{S}^3$ are precisely the Hopf vector fields, and no others.
\end{theorem}
On the other hand, Reznikov compared this functional to the topology of an Euclidean hypersurface. Let $M^{n+1}$ be a smooth closed oriented immersed hypersurface in $\mathbb{R}^{n+2}$, endowed with the induced metric, and let $\mathcal{S}=\sup_{x\in M}\|S_x\|=\sup_{x\in M}|\lambda_i(x)|$, where $S_x$ is the second fundamental operator in $T_xM$, and $\lambda_i(x)$ are the principal curvatures.
\begin{theorem}[Reznikov, \cite{reznikov}]
\label{rez}
For any unit vector field $\vec{v}$ on $M$ we have
$$\vol(\vec{v})-\vol(M)\geq
\frac{\vol(\mathbb{S}^{n+1})}{\mathcal{S}}|\deg(\nu)|,$$
where $\deg(\nu)$ is the degree of the Gauss map $\nu:M\rightarrow\mathbb{S}^{n+1}$.
\end{theorem}

In this short note, we take an odd dimensional hypersurface $M^{2n+1}$ and relate the energy of a unit vector field $\vec{v}$ to the topology of the immersion of $M$, by means of the degree of the Gauss map. Our main theorem reads
\begin{theoremA}\label{teo-energy} For a unit vector field on a closed oriented Euclidean hypersurface $M^{2n+1}$, 
\begin{equation}
E(\vec{v})\geq  C(n)\frac{ |\!\deg(\nu)| \vol(\mathbb{S}^{2n+1})}{S^{[2n-1]}} + \frac{2n+1}{2}\vol(M^{2n+1})\nn
\end{equation}
where $S^{[2n-1]}$ and $C(n)$ are constants depending on the immersion of $M$ and on $n$ (their precise definition will be given later). 
\end{theoremA}

Theorem A provides a topological obstruction to small values of the energy in a Riemannian manifold, specifically in a hypersurface in the Euclidean space. Two non-homotopic immersions will possess two different normal degrees; the bigger this value, the bigger the energy of a given unit vector field. As far as the authors know, this is the first connection between the topology of an immersion and the energy of a unit vector field. 

A special case is the unit sphere $\mathbb{S}^{2n+1}$.  Borreli {\it et al} \cite{borreli} constructed a family of unit vector fields on $\mathbb{S}^{2n+1}$ with energy converging to the energy of a radial vector field. Its value turned out to be the infimum for the energy of unit vector fields without singularities.

\begin{theorem}[Borreli, Brito and Gil-Medrano, \cite{borreli}]\label{thm-borreli-brito-olga} The infimum of $E$ among all globally defined unit smooth vector fields
of the sphere $\mathbb{S}^{2n+1}$ $(n\geq 2)$ is
\begin{equation}
\left(\frac{2n+1}{2} + \frac{n}{2n-1}  \right) \vol(\mathbb{S}^{2n+1}).
\end{equation}
This value is not attained by any globally defined unit smooth vector field.
\end{theorem}
 
By a theorem of Hopf, $\deg(\nu)$ is the same for homotopic immersions of a compact hypersurface. Thus, it is interesting to look at spheres of different radii. 
Applying the procedure from theorem A to a sphere of radius $r$, we have

\begin{corollary}
\label{coro}
For $r>0$, let $\mathbb{S}^{2n+1}(r)$ be immersed in $\mathbb{R}^{2n+2}$ with normal degree one. Then 
\begin{equation}
E(\vec{v})\geq \left(\frac{2n+1}{2}r^{2n+1} + \frac{n}{2n-1} r^{2n-1}  \right) \vol(\mathbb{S}^{2n+1}).\nn
\end{equation}
\end{corollary}
Consequently, when $r=1$ we recover the value from theorem \ref{thm-borreli-brito-olga}. 
This inequality is well known from the literature, see \cite{brito-pawel} for a further discussion on the energy of vector fields possessing isolated singularities, as well as a proof of a general inequality regarding ${\rm Ricci}(\vec{v},\vec{v})$ on a given compact Riemanninan manifold.

We also discuss a list of functionals having properties similar to the total bending of flows, and determine a lower value for each one of them depending again on $\deg(\nu)$. Let $\vec{v}$ be a unit vector field on a compact oriented Riemannian manifold $M^m$. For every $1\leq k\leq m-1$, define 
\begin{equation}
\label{hotb}
\mathcal{B}_k(\vec{v}) = \int_M \underbrace{\|\nabla \vec{v}\wedge \cdots\wedge \nabla \vec{v} \|^2}_{\text{$k$-times}}.
\end{equation}
If $\sigma_{2n}$ denotes the $2n$-th elementary symmetric function, and $\mathcal{V}$ is the restriction of $\nabla\vec{v}$ to $V^{\perp}$ then our last theorem reads
\begin{theoremB}
Let $M^{2n+1}$ be a compact oriented Riemannian manifold, and let $\vec{v}$ be a unit vector field on $M$. Then
\begin{equation}
\label{lower-b-Bn}
\mathcal{B}_n(\vec{v})\geq {2n \choose n}\int_M |\sigma_{2n}(\mathcal{V}) |.
\end{equation}
Furthermore, when $M^{2n+1}$ is a closed Euclidean hypersurface, 
\begin{equation}
\label{lower-b-Bn-hypers}
\mathcal{B}_n(\vec{v})\geq \frac{|\deg(\nu)|}{\mathcal{S}}{2n \choose n} \vol(\mathbb{S}^{2n+1}), 
\end{equation}
where $\mathcal{S}$ is the aforementioned constant.
\end{theoremB}

As a consequence, we deduce the following
\begin{corollary}
\label{hopf4n1}
Hopf vector fields minimize $\mathcal{B}_n$ on $S^{2n+1}$.
\end{corollary}

Hopf vector fields on $\mathbb{S}^3$ are absolute minima of the energy, \cite{brito}, and on $\mathbb{S}^{2n+1}$, for $n\geq 2$, they are critical but unstable points of the energy functional, \cite{wood}. Despite these properties on higher dimensional spheres, the functionals $\mathcal{B}_n$ on $S^{2n+1}$ are an attempt to provide a list of functionals that are minimized by Hopf vector fields, having similar features when compared to the energy and/or total bending. This result should also be compared to a mean curvature correction of the total bending provided by the first author on \cite{brito}. 

\section{Curvature integrals for a closed hypersurface}
\label{eta}

The proofs of theorem A, corollary \ref{coro} and of \ref{lower-b-Bn-hypers} rely on a list of curvature integrals, described in this Section. 

We may assume that $M^{2n+1}$ is oriented, so the normal map $\nu:M^{2n+1}\rightarrow \mathbb{S}^{2n+1}$,  $\nu(x)=N(x)$, is well defined, where $N$ is a unitary normal field. Let $\langle\cdot,\cdot\rangle$ be the induced Riemannian metric on $M$. Let $\vec{v}:M\rightarrow TM$ be a smooth unit vector field on $M$, and take the orthonormal basis $\{e_1,\dots,e_{2n}, e_{2n+1}:=\vec{v}\}$ at each point $x\in M$. We fix some notation: for $1\leq A,B\leq2n+1$, set $h_{AB}=\langle S(e_A),e_B\rangle$ and $a_{AB}=\langle\nabla_{e_B}\vec{v},e_A\rangle$; it follows that  $a_{2n+1\, B} =\langle\nabla_{e_B}\vec{v},\vec{v}\rangle = 0$, for all $B$. For a real number $t>0$, define $\varphi^{\vec{v}}_{t}:  M^{2n+1}\to \mathbb{S}^{2n+1}(\sqrt{1+t^2})$, by $\varphi^{\vec{v}}_{t}(x) = \nu(x) + t\vec{v}(x)$. With respect to the aforementioned basis of $T_xM$ and setting $\left\lbrace e_1, \dots, e_{2n}, \frac{\vec{v}}{\sqrt{1+t^2}} - t\frac{N}{\sqrt{1+t^2}}\right\rbrace$ as an orthonormal basis of $T_{\varphi^{\vec{v}}_{t}(x)} \mathbb{S}^{2n+1}(\sqrt{1+t^2})$, we have that
\[
d\varphi^{\vec{v}}_{t} = 
  \left(\begin{array}{ccc|c}
    \multicolumn{3}{c|}{\multirow{3}{*}{\raisebox{2pt}{$h_{ij} + ta_{ij}$}}}     & h_{2n+1\;1} + ta_{1\, 2n+1}           \\ 
    &       &            & {\vdots}    \\
    &       &            & h_{2n+1\;2n} + ta_{2n\, 2n+1}          \\ \hline
    \sqrt{1+t^2}h_{2n+1\;1}      & \cdots & \sqrt{1+t^2}h_{2n+1\; 2n}  & \sqrt{1+t^2}h_{2n+1\; 2n+1}
  \end{array}\right).
\]
Multilinearity of determinant simplifies computations concerning an explicit formula for $\det(d\varphi^{\vec{v}}_{t})$ written in terms of the second fundamental form of $M$ and components depending on the normal bundle of $\vec{v}$. 
\begin{eqnarray*}
	\det((h_{AB})+t(a_{AB})) & = & \sum_{\sigma}\mbox{sgn}(\sigma)\Big(h_{1\,\sigma(1)}\cdots h_{2n+1\,\sigma(2n+1)}	 \\ \nonumber
 	& + & t\sum_i h_{1\,\sigma(1)}\cdots a_{i\,\sigma(i)}\cdots h_{2n+1\,\sigma(2n+1)}	\\ \nonumber
 	& + &  t^2\sum_{i<j} h_{1\,\sigma(1)}\cdots a_{i\,\sigma(i)}\cdots a_{j\,\sigma(j)}\cdots h_{2n+1\,\sigma(2n+1)} \\	\nonumber
	&\vdots&\nn\\
	&+& t^{2n}a_{1\,\sigma(1)} \cdots a_{2n\,\sigma(2n)} h_{2n+1\,\sigma(2n+1)}\Big)	\nonumber
\end{eqnarray*}
Therefore,
\begin{equation}
\det(d\varphi^{\vec{v}}_{t}) = \sqrt{1+t^2} \sum_{k=0}^{2n}\eta_k t^k \nn
\end{equation}
where
\begin{eqnarray}
\eta_0 &=& \det(h_{AB})\nonumber\\
\eta_1 &=& \sum_{\sigma}\mbox{sgn}(\sigma) \sum_i h_{1\,\sigma(1)}\cdots h_{i-1\,\sigma(i-1)}a_{i\,\sigma(i)}h_{i+1\,\sigma(i+1)}\cdots h_{2n+1\,\sigma(2n+1)} \nonumber\\
\eta_2 &=& \sum_{\sigma}\mbox{sgn}(\sigma) \sum_{i<j} h_{1\,\sigma(1)}\cdots h_{i-1\,\sigma(i-1)}a_{i\,\sigma(i)}h_{i+1\,\sigma(i+1)}\cdots\nonumber\\
&\cdots& h_{j-1\,\sigma(j-1)}a_{j\,\sigma(j)}h_{j+1\,\sigma(j+1)}\cdots h_{2n+1\,\sigma(2n+1)} \nonumber\\
&\vdots & \nonumber\\
\eta_{2n} &=&	\sum_{\sigma}\mbox{sgn}(\sigma)a_{1\,\sigma(1)} \cdots a_{2n\,\sigma(2n)} h_{2n+1\,\sigma(2n+1)}. \nonumber
\end{eqnarray}

The fact that $\det((h_{AB}) + t (a_{AB})) = \sum_k \eta_k t^k$ shows that, for $1\leq k \leq 2n$, $\eta_k$ does not depend on the choice of basis.

Since $\varphi^{\vec{v}}_{t}$ and $\nu$ are homotopic, we have that $\deg(\varphi^{\vec{v}}_{t}) = \deg(\nu)$. On the other hand, if $f:X\to Y$ is a smooth map between two manifolds of the same dimension, say $m$, then for any $m$-form $\omega$ on $Y$, the degree formula reads $\int_Xf^*(\omega) = \deg(f)\int_Y\omega$. By change of variables we conclude that  
\begin{equation}
\label{integral-formula}
\int_{M} \eta_k=
\begin{cases}
\deg(\nu){n \choose k/2} {\rm vol}(\mathbb{S}^{2n+1}),
  &\mbox{if}\; k \; \mbox{is}\;\mbox{even},\\
0, & \mbox{if}\; k \;\mbox{is}\;\mbox{odd}.
\end{cases}
\end{equation}

These invariants have been computed and applied to questions concerning geometry of foliations on hypersurfaces, \cite{icaro}. 

\section{Proof of theorem A}
We start our approach by defining a list of numbers regarding wedge products of the shape operator of $M$ and their restriction to a list of vectors on a point. 
\begin{definition}
\label{norm-max-multiv}
If $\{u_1,\dots,u_{2n+1}\}$ is an orthonormal basis at $x\in M$, then, for each $1\leq A\leq2n+1$,
\begin{equation}
S^{[A]} = \sup_{1\leq i_1,\dots, i_{A}\leq 2n+1;\; x\in M}\{\|S(u_{i_1})\wedge \cdots\wedge S(u_{i_A})\|_{\infty}\},\nn
\end{equation}
where $\|\cdot\|_{\infty}$ denotes the maximum norm,  naturally extended to $\Lambda^A(M)$. 
\end{definition}
We notice that $S^{[A]}$ is well defined since $M$ is a compact hypersurface. With respect to $\{e_1,\dots,e_{2n}, \vec{v}\}$, the energy functional is written as 
\begin{eqnarray}
E(\vec{v})&=&\frac{1}{2}\int_M\|\nabla\vec{v}\|^2+\frac{2n+1}{2}\vol(M)\nn\\
&=&\frac{1}{2}\int_M\left(\sum_{A,B}a_{AB}^2\right)+\frac{2n+1}{2}\vol(M),\nn
\end{eqnarray}
this means that $\eta_2$ is the natural choice among all $\eta_k$ in order to determine a lower value for $E(\vec{v})$. From the last Section, we have
\begin{eqnarray}
\int_M\eta_2 &=& \sum_{\sigma}\mbox{sgn}(\sigma) \sum_{i<j} h_{1\,\sigma(1)}\cdots h_{i-1\,\sigma(i-1)}a_{i\,\sigma(i)}h_{i+1\,\sigma(i+1)}\cdots\nonumber\\
&\cdots& h_{j-1\,\sigma(j-1)}a_{j\,\sigma(j)}h_{j+1\,\sigma(j+1)}\cdots h_{2n+1\,\sigma(2n+1)} \nonumber\\
&=& n\deg(\nu) \vol(S^{2n+1}) \nn
\end{eqnarray}

By the definition of the matrix $(a_{AB})$ and the fact that two of its entries appear in the above summation, we are able to display all $2\times 2$ minors of $(a_{AB})$, times a minor depending on the second fundamental form matrix of $M$. Precisely,   
\begin{equation}
\label{eta2}
\int_M\eta_2 = \int_M\left(\sum_{i,j,k,l}\left|\begin{array}{cc}
a_{ij} & a_{ik} \\ 
a_{lj} & a_{lk}
\end{array}\right| \det(\widehat{H}^{ij}_{lk}) + \sum_{i,j,l}\left|\begin{array}{cc}
a_{ij} & a_{j\, 2n+1} \\ 
a_{lj} & a_{l\, 2n+1}
\end{array}\right| \det(\widehat{H}^{i\, 2n+1}_{l\, 2n+1})\right),
\end{equation}
where $\widehat{H}^{AB}_{CD}$ is a $(2n-1)\times (2n-1)$ matrix which comes from $(h_{AB})$ by removing the $A$-th and $C$-th lines, and $B$-th and $D$-th columns. 

For $0\leq k \leq 2n$ the functions $\eta_k$ are invariant under any change of basis. We may assume that $\eta_2$ is computed with respect to a basis that diagonalizes the second fundamental form matrix of $M$. In this new basis, we have a matrix $(\widetilde{a}_{AB})$ and its last line might be different from zero. So we write 
\begin{eqnarray}
\int_M\eta_2 &=& \int_M \sum_{A<D, B<C}\left|\begin{array}{cc}
\widetilde{a}_{AB} & \widetilde{a}_{AC} \\ 
\widetilde{a}_{DB} & \widetilde{a}_{DC}
\end{array}\right| \det(\widehat{H}^{AB}_{CD})\nonumber\\
&=& \int_M \sum_{A<D, B<C}\left|\begin{array}{cc}
\widetilde{a}_{AB} & \widetilde{a}_{AC} \\ 
\widetilde{a}_{DB} & \widetilde{a}_{DC}
\end{array}\right| \prod_{F\neq A, C}h_{FF}.\nonumber
\end{eqnarray}
From our definition \ref{norm-max-multiv}, $\prod_{F\neq A, C}h_{FF}\leq S^{[2n-1]}$. Employing the inequality $\left|\begin{array}{cc}
a & b \\ 
c & d
\end{array}\right| \leq \frac{1}{2} (a^2+b^2+c^2+d^2)$, we have  
\begin{equation}
\int_M\eta_2 \leq \frac{1}{2} \int_M \sum_{A<D, B<C} \left(\widetilde{a}_{AB}^2 + \widetilde{a}_{AC}^2 +
\widetilde{a}_{DB}^2 + \widetilde{a}_{DC}^2\right) S^{[2n-1]}.
\nonumber
\end{equation}

Now comes the crucial distinction between theorem A and corollary \ref{coro}; that is, between the constant $C(n)$ obtained for an arbitrary hypersurface and by restricting to a sphere of radius $r$. For an arbitrary closed hypersurface $M^{2n+1}$, we are able to diagonalize $(h_{AB})$ but we can not control which entries in $(\widetilde{a}_{AB})$ will remain different from zero. This implies that when we count its squared entries in the last inequality we end up with $2n$ elements. Hence
\begin{equation}
\int_M\eta_2 \leq n S^{[2n-1]} \int_M \sum_{A,B}\widetilde{a}_{AB}^2 \leq n S^{[2n-1]} \int_M \|\nabla\vec{v}\|^2.
\nonumber
\end{equation}

On the other hand, if $M=\mathbb{S}^{2n+1}(r)$, then $(h_{AB})$ is $1/r$ times the identity matrix, $(h_{AB}) = (1/r)\cdot I$. This means that for any choice of basis $(h_{AB})$ is a diagonal matrix. In particular, we can employ the orthonormal basis $\{e_1,\dots,e_{2n}, e_{2n+1}:=\vec{v}\}$, and $(a_{AB})$ is the same as before, i.e. $(a_{2n+1\ B}) = \langle \nabla_{e_B}\vec{v}, \vec{v} \rangle = 0$ for all $1\leq B \leq 2n+1$. In this case, when we count the terms in $\sum_{A<D, B<C} \left(a_{AB}^2 + a_{AC}^2 +
a_{DB}^2 + a_{DC}^2\right)$ we have a smaller number, $2n-1$. Thus 
\begin{equation}
\int_{\mathbb{S}^{2n+1}(r)}\eta_2  \leq \frac{2n-1}{2} S^{[2n-1]} \int_{\mathbb{S}^{2n+1}(r)} \|\nabla\vec{v}\|^2.
\nonumber
\end{equation}

Therefore,
\begin{equation}
E(\vec{v})\geq  C(n)\frac{ |\!\deg(\nu)| \vol(\mathbb{S}^{2n+1})}{S^{[2n-1]}} + \frac{2n+1}{2}\vol(M^{2n+1}),\nn
\end{equation}
where
\begin{equation}
C(n) =
\begin{cases}
\frac{n}{2n-1},
  &\text{if}\ M^{2n+1} = \mathbb{S}^{2n+1}(r),
  \\
\frac{1}{2}, & \text{otherwise}.
\end{cases}\nn
\end{equation}

Notice that when $ M^{2n+1} = \mathbb{S}^{2n+1}(r)$, its volume is $r^{2n+1}$ times the volume of $\mathbb{S}^{2n+1}$. In addition, $S^{[2n-1]} = r^{-2n+1}$, because $(h_{AB}) = (1/r)\cdot I$. This concludes the proof of corollary \ref{coro}.

\section{Known results and a list of new functionals}
This section is intended to show some direct consequences of equation \ref{integral-formula}, as well as proving theorem B. We start by considering the integral of the last invariant $\eta_{2n}$ and then we obtain the main ingredient in the proof of theorem \ref{rez} from \cite{reznikov}. Finally, we discuss the list of ``higher order" total bending functionals on Riemannian manifolds, and study their lower bounds on closed Euclidean hypersurfaces. 

\subsection{Another view on theorem \ref{rez}}

By definition, the last line of $(a_{AB})$ is zero ($\vec{v}$ is unitary), so we let $m_A$, $1\leq A\leq 2n+1$, denote its $2n$-minors. Thus, $\eta_{2n} = \sum_A (-1)^{A+1}h_{2n+1 A} m_A$, which implies 
\begin{eqnarray}
|\eta_{2n}| &\leq& \left(\sum_A h_{2n+1 A}^2\right)^{1/2}\left(\sum_A m_{A}^2\right)^{1/2} = \|S(\vec{v})\| \left(\sum_A m_{A}^2\right)^{1/2}\nn\\
&\leq& \mathcal{S}\left(\sum_A m_{A}^2\right)^{1/2}\nn\\
&\leq& \mathcal{S} \left(\sqrt{\det(I+(\nabla\vec{v})(\nabla\vec{v})^t)} - 1\right).\nn
\end{eqnarray}
For the last inequality, see for example Lemma 1 in \cite{reznikov}. The proof of the theorem \ref{rez} is finished by integration over $M$. 
\subsection{Higher order functionals: proof of theorem B}
The fact that $\vec{v}$ is unitary implies that $\mathcal{B}_m(\vec{v}) = 0$, and $\mathcal{B}_1(\vec{v})$ is, up to a constant, the total bending of $\vec{v}$. All $\mathcal{B}$ functionals can be written as integrals of functions of $k$-minors from the matrix $(a_{AB})$, 
\begin{equation}
\mathcal{B}_k(\vec{v}) = \int_M \sum_{\substack{A_1<\cdots< A_k \\ B_1<\cdots< B_k}} {\det}^2\left(a^{A_1\cdots A_k}_{B_1\cdots B_k} \right).\nn
\end{equation}

Assume that $M$ has dimension $2n+1$. When $1\leq i,j\leq 2n$, the matrix $(a_{ij}) = (\langle \nabla_{e_j} \vec{v}, e_i \rangle)$ describes the behavior of the distribution normal to $\vec{v}$. We are going to compare its determinant to the integrand in the functional $\mathcal{B}_n$. 

If we omit the $n$-minors having at least one element of the type $a_{i\ 2n+1} = \langle \nabla_{\vec{v}} \vec{v}, e_i \rangle$, then
\begin{equation}
\mathcal{B}_n(\vec{v}) = \int_M \sum_{\substack{A_1<\cdots< A_n \\ B_1<\cdots< B_n}} {\det}^2\left(a^{A_1\cdots A_n}_{B_1\cdots B_n} \right) \geq \int_M \sum_{\substack{i_1<\cdots< i_n \\ j_1<\cdots< j_n}} {\det}^2\left(a^{i_1\cdots i_n}_{j_1\cdots j_n} \right). \nn
\end{equation}

We follow Sections 3 and 4 from \cite{BCN}, and Chapter IV of \cite{tese-pc}. In general, the distribution normal to $\vec{v}$ is not integrable. Even though $(a_{ij})$ is not symmetric, we can find a local basis in which $(a_{ij})$ has the form of a upper triangular matrix, 
\[(a_{ij}) = 
    \left(
    \begin{array}{cccccccccc}
    \lambda_1  &\ast & \cdots &&&&&&\cdots& \ast \\
    0  & \ddots    &   & &&&&&&\vdots\\
    \vdots  &  &\ \lambda_r & \ast & \cdots  \\
      &  & 0   &             \\
      & & \vdots & & x_1 & -y_1 &\ast & \cdots \\
      & & & & y_1&\ x_1\\
      & & & & 0& & \ddots &&&\vdots\\
      & & & & \vdots & & & \ddots & & \ast \\
     \vdots & & & & & & & & x_s & -y_s\\
     0 &\cdots & & & & & \cdots& 0&y_s &\ x_s
    \end{array}
    \right).
\]  

Let $D = \text{diag}(|\lambda_1|,\dots,|\lambda_r|,\sqrt{x_1^2 + y_1^2}, \sqrt{x_1^2 + y_1^2}, \dots, \sqrt{x_s^2 + y_s^2}, \sqrt{x_s^2 + y_s^2})$ be a $2n\times 2n$ diagonal matrix. By construction, $\det D \geq |\det (a_{ij})|$. On the other hand, 
$$
\sum_{\substack{i_1<\cdots< i_n \\ j_1<\cdots< j_n}} {\det}^2\left(a^{i_1\cdots i_n}_{j_1\cdots j_n} \right) \geq \sum_{\substack{i_1<\cdots< i_n \\ j_1<\cdots< j_n}} {\det}^2\left(D^{i_1\cdots i_n}_{j_1\cdots j_n} \right),
$$
simply because most elements above the main diagonal in $(a_{ij})$ are different from zero, which is the case for $D$. 

The main result of Section 3 in \cite{BCN} is the ``Fundamental Lemma", which states an inequality between the volume of a $2n\times 2n$ diagonal matrix and the sum of its even elementary symmetric functions. In proving this lemma, the authors deduced the following inequality (cf. second inequality on page 307 of \cite{BCN}; we also refer to \cite{tese-pc}, equations (IV.16) and (IV.21) pages 55 and 56, respectively) 
$$
\sum_{\substack{i_1<\cdots< i_n \\ j_1<\cdots< j_n}} {\det}^2\left(D^{i_1\cdots i_n}_{j_1\cdots j_n} \right) \geq {2n \choose n} \det D.
$$
Since $\sigma_{2n}(\mathcal{V}) = \det (a_{ij})$, this completes the proof of \ref{lower-b-Bn}. 

The simplest invariant depending on $\vec{v}$ in Section \ref{eta} is $\eta_{2n}$. If $\mathcal{S}$ is the number defined in theorem \ref{rez}, then both of them combines to establish \ref{lower-b-Bn-hypers}.

\subsection{Hopf flows on $\mathbb{S}^{2n+1}$}
We proceed to prove Corollary \ref{hopf4n1}. Hopf vector fields are tangent to the fibers of $\mathbb{S}^1\hookrightarrow \mathbb{S}^{2n+1}\to \mathbb{C}{\rm P}^n$, and the matrix associated to the second fundamental form of a given $H$ is a $2n\times 2n$ matrix of the type
\[
(a_{ij})=\left(
\begin{array}{ccccc}
 \frame{$\begin{array}{cc} 0 & -1 \\ 1 & 0 \\ \end{array}$} &  &   \\
   & \ddots & \\
 &  & \frame{$\begin{array}{cc} 0 & -1 \\ 1 & 0 \\ \end{array}$} 
\end{array}\right),
\]    
having $n$ blocks of the form $\left(\begin{array}{cc} 0 & -1 \\ 1 & 0 \\ \end{array}\right)$ in a diagonal and zeros everywhere else. In computing $\mathcal{B}_n(H)$ on $\mathbb{S}^{2n+1}$, we need to count how many $n$-minors different from zero the above matrix has. A given $n\times n$ sub-matrix is obtained by removing $n$ rows and $n$ columns of $(a_{ij})$. At the end of this process, all rows and columns of this $n\times n$ sub-matrix have exactly one element, $\pm 1$. By rearranging the rows in this sub-matrix we get a $n\times n$ diagonal matrix, and this matrix has determinant $\pm 1$. Thus, we have ${2n \choose n}$ non-zero $n$-minors, each one evaluating $\pm 1$. Therefore, $\mathcal{B}_n(H) = {2n \choose n}\vol(\mathbb{S}^{2n+1})$. 

\section*{Acknowledgment}
We thank an anonymous referee for valuable suggestions and a critical look at a previous version of this manuscript.

\end{document}